\documentclass[12pt,twoside,a4paper]{article}
\usepackage[english]{babel}
\usepackage{amssymb}
\usepackage[T1]{fontenc}
\usepackage{amsmath,epsf}
\newtheorem{Th}{Theorem}
\newtheorem{lemma}{Lemma}

\newtheorem{Prop}{Proposition}

\newcommand{\E}{\mathbb{E}}
\newcommand{\F}{\mathcal{F}}
\newcommand{\R}{\mathbb{R}}
\newcommand{\Z}{\mathbb{Z}}
\newcommand{\N}{\mathbb{N}}

\newcommand{\U}{\mathcal{U}}
\newcommand{\V}{\mathbb{V}}

\renewcommand{\P}{\mathbb{P}}

\newcounter{tictac}

\def\1{\,\rlap{\mbox{\small\rm 1}}\kern.15em 1}
\def\ind#1{\1_{#1}}
\def\build#1_#2^#3{\mathrel{\mathop{\kern 0pt#1}\limits_{#2}^{#3}}}
\def\tend#1#2{\build\hbox to 12mm{\rightarrowfill}_{#1\rightarrow #2}^{a.s.}}

\def\converge#1#2#3{\build\hbox to
15mm{\rightarrowfill}_{#1\rightarrow #2}^{\hbox{\scriptsize #3}}}
\begin{document}

\begin{center}
{\large
    {\sc
On the asymptotic normality of frequency polygons for strongly mixing spatial processes
    }
}
\bigskip

 Mohamed EL MACHKOURI

\medskip
\scriptsize{Laboratoire de Math\'ematiques Rapha\"el Salem\\
UMR CNRS 6085, Universit\'e de Rouen (France)\\
\emph{mohamed.elmachkouri@univ-rouen.fr}}
\end{center}

{\renewcommand\abstractname{Abstract}
\begin{abstract}
\baselineskip=18pt This paper establishes the asymptotic normality of frequency polygons in the context of stationary strongly mixing random fields indexed by $\Z^d$. Our method allows us to consider only minimal conditions on the width bins and provides a simple criterion on the mixing coefficients. In 
particular, we improve in several directions a previous result by Carbon, Francq and Tran (2010).\\
\\
{\em AMS Subject Classifications} (2000): 62G05, 62G07, 60G60.\\
{\em Key words and phrases:} Central limit theorem, spatial processes, random fields, nonparametric density estimation, frequency polygon, histogram, mixing.\\
\end{abstract}

\thispagestyle{empty}
\baselineskip=18pt
\section{Introduction and notations}
The frequency polygon is a density estimator based on the histogram. It has the advantage to be conceptually and computationaly simple since it just consists of 
linking the mid points of the histogram bars but its simplicity is not the only interest. In fact, for time series, 
Scott \cite{Scott1985} as shown that the rate of convergence of frequency polygon is superior to the histogram for smooth densities. 
For other references on non-spatial density estimation based on the frequency polygon, one can refer to Beirlant et al. \cite{Beirlant-Berlinet-Gyorfi1999} 
and Carbon et al. \cite{Carbon-Garel-Tran1997}. To our knowledge, the only references in the spatial context are Carbon \cite{Carbon2006} and 
Carbon et al. \cite{Carbon-Francq-Tran2010} for strongly mixing random fields indexed by $\Z^d$ and Bensaid and Dabo-Niang \cite{Bensaid-DaboNiang2010} 
for strongly mixing random fields indexed by $\R^d$. In \cite{Carbon-Francq-Tran2010} the asymptotic normality of the frequency polygon estimator is 
obtained under interleaved conditions on the width bin and the strong mixing coefficients. In this paper, we provide a simple criterion on the mixing coefficients 
for the frequency polygon to satisfy the central limit theorem when only minimal conditions on the width bin (see Assumption $\textbf{(A2)}$ below) are assumed. Our main result (Theorem 
\ref{convergence-loi}) improve Theorem 4.1 in \cite{Carbon-Francq-Tran2010} in several directions. Our approach which is based on the Lindeberg's method seems to be superior to the so-called Bernstein's blocking method used 
in \cite{Carbon-Francq-Tran2010} but also in many others papers on nonparametric estimation for random fields 
(see \cite{Bensaid-DaboNiang2010}, \cite{Carbon-Hallin-Tran}, \cite{Carbon-Tran-Wu}, \cite{Hallin-Lu-Tran-2001}, \cite{Tran1990}). For another application of 
the Lindeberg's method, one can refer to El Machkouri \cite{elmachkouri2011} where the asymptotic normality of the Parzen-Rosenblatt kernel density estimator is proved for 
strongly mixing random fields improving a previous result by Tran \cite{Tran1990} obtained also by the Bernstein's blocking method an coupling arguments. Note that 
the central limit theorem in \cite{elmachkouri2011} is obtained for random fields observed on squares $\Lambda_n$ of $\Z^d$ but actually the result still holds if the regions 
$\Lambda_n$ are only assumed to have cardinality going to infinity as $n$ goes to infinity. In particular, it is not neccessary to impose any condition on the boundary of $\Lambda_n$. \\ 
Let $d$ be a positive integer and let $(X_i)_{i\in\Z^d}$ be a stationary real random field defined on some probability space $(\Omega, \F,\P)$ with an unknown 
marginal density $f$. For any finite subset $B$ of $\Z^d$, denote $\vert B\vert$ the number of elements in $B$. 
In the sequel, we assume that we observe $(X_i)_{i\in\Z^d}$ on a sequence $(\Lambda_n)_{n\geq 1}$ of finite subsets of $\Z^d$ such that 
$\vert\Lambda_n\vert$ goes to infinity as $n$ goes to infinity. We lay emphasis on that we do not impose any condition on the boundary of the regions $\Lambda_n$.
Given two $\sigma$-algebras $\U$ and $\mathcal{V}$, we recall the $\alpha$-mixing coefficient introduced by Rosenblatt \cite{Ros} and defined by 
$$
\alpha(\U,\mathcal{V})=\sup\{\vert\P(A\cap B)-\P(A)\P(B)\vert\, ,\,A\in\U,\,B\in\mathcal{V}\}
$$ 
and the $\rho$-mixing coefficient introduced by Kolmogorov and Rozanov \cite{Kolmogorov--Rozanov1960} and defined by
$$
\rho(\U,\mathcal{V})=\sup\left\{\frac{\vert\textrm{Cov}(f,g)\vert}{\|f\|_2\|g\|_2}\, ,\,f\in\mathbb{L}^2(\mathcal{U}),\,g\in\mathbb{L}^2(\mathcal{V})\right\}.
$$
It is well known that $4\alpha(\U,\mathcal{V})\leq\rho(\U,\mathcal{V})$ (see \cite{Bradley2005}). For any $\tau$ in $\N^{\ast}\cup\{\infty\}$ and any positive integer $n$, we consider the mixing coefficients $\alpha_{1,\tau}(n)$ and  $\rho_{1,\tau}(n)$ defined by
\begin{align*}
\alpha_{1,\tau}(n)&=\sup\,\{\alpha(\sigma(X_k),\F_{B}),\,k\in\Z^d,\,\vert B\vert\leq\tau,\,\Xi(B,\{0\})\geq n\},\\
\rho_{1,\tau}(n)&=\sup\,\{\rho(\sigma(X_k),\F_{B}),\,k\in\Z^d,\,\vert B\vert\leq\tau,\,\Xi(B,\{0\})\geq n\}
\end{align*}
where $\F_{B}=\sigma(X_i\,;\,i\in B)$ for any subset $B$ of $\Z^d$ and the distance $\Xi$ is defined for any subsets $B_1$ and $B_2$ of $\Z^d$ by 
$\Xi(B_{1},B_{2})=\min\{\vert i-j\vert,\,i\in B_{1},\,j\in B_{2}\}$ where $\vert i-j\vert=\max_{1\leq k\leq d}\vert i_k-j_k\vert$ for any $i=(i_1,...,i_d)$ and 
$j=(j_1,...,j_d)$ in $\Z^d$. We say that the random field $(X_i)_{i\in\Z^d}$ is $\alpha$-mixing or $\rho$-mixing if $\lim_{n\to\infty}\alpha_{1,\tau}(n)=0$ or $\lim_{n\to\infty}\rho_{1,\tau}(n)=0$ for some $\tau$ in $\N^{\ast}\cup\{\infty\}$ respectively.\\
Let $(b_n)_{n\geq 1}$ be a sequence of positive real numbers going to zero when $n$ goes to infinity. For each $n$ in $\N^{\ast}$, we consider the partition 
$(\textrm{I}_{n,k})_{k\in\Z}$ of the real line defined for any $k$ in $\Z$ by $\textrm{I}_{n,k}=[(k-1)b_n,kb_n)$. Let $(n,k)$ be fixed in $\N^{\ast}\times\Z$ and let 
$\textrm{I}_{n,k}$ and $\textrm{I}_{n,k+1}$ be two adjacent histogram bins. Denote $\nu_{n,k}$ and $\nu_{n,k+1}$ the numbers of observations falling in these intervals 
respectively. The values of the histogram in these previous bins are given by $\nu_{n,k}(\vert\Lambda_n\vert b_n)^{-1}$ and $\nu_{n,k+1}(\vert\Lambda_n\vert b_n)^{-1}$ 
and the frequency polygon $f_{n,k}$ is defined for $x\in \textrm{J}_{n,k}:=\left[\left(k-\frac{1}{2}\right)b_n,\left(k+\frac{1}{2}\right)b_n\right)$ by
$$
f_{n,k}(x)=\left(\frac{1}{2}+k-\frac{x}{b_n}\right)\frac{\nu_{n,k}}{\vert\Lambda_n\vert b_n}+\left(\frac{1}{2}-k+\frac{x}{b_n}\right)\frac{\nu_{n,k+1}}{\vert\Lambda_n\vert b_n}.
$$
Define $Y_{i,s}=\ind{X_i\in I_{n,s}}$ for any $s$ in $\Z$, then 
$$
f_{n,k}(x)=\frac{1}{\vert\Lambda_n\vert b_n}\sum_{i\in\Lambda_n}a_k(x)Y_{i,k}+\overline{a}_k(x)Y_{i,k+1}
$$
where  $a_s(u)=\frac{1}{2}+s-\frac{u}{b_n}$ and $\overline{a}_{s}(u)=a_{-s}(-u)$ for any $s$ in $\Z$ and any $u$ in $\textrm{J}_{n,s}$. Finally, we consider 
the normalized frequency polygon estimator $f_n$ defined for any $x$ in $\R$ such that $f(x)>0$ by
$$
f_n(x)=\sum_{k\in\Z}\frac{f_{n,k}(x)}{\sigma_{n,k}(x)}\ind{\textrm{J}_{n,k}}(x)\quad\textrm{where}\quad\sigma^2_{n,k}(x)=\left(\frac{1}{2}+2\left(k-\frac{x}{b_n}\right)^2\right)f(x).
$$
Our main results are stated in Section $2$ and the proofs are given in Section $3$.
\section{Main results}
We consider the following assumptions.
\begin{itemize}
\item[\textbf{(A1)}] The density function $f$ is differentiable and its derivative $f^{'}$ is locally bounded.
\item[\textbf{(A2)}] $b_n$ goes to zero such that $\vert\Lambda_n\vert b_n$ goes to infinity as $n$ goes to infinity.
\item[\textbf{(A3)}] $\sup_{j\neq 0}\P\left(X_0\in I_{n,s},X_j\in I_{n,t}\right)=O(b_n^2)$ for any $s$ and $t$ in $\Z$.
\item[\textbf{(B3)}] $\P\left(X_0\in I_{n,s},X_j\in I_{n,t}\right)=o(b_n)$ for any $s$, $t$ and $j$ in $\Z$.
\end{itemize}
\textbf{Remark 1}. Obviously $\textbf{(B3)}$ is weaker than $\textbf{(A3)}$. Moreover, if the joint density $f_{0,j}$ of $(X_0,X_j)$ exists then $\textbf{(A3)}$ is true by assuming that $\sup_{j\neq 0} f_{0,j}$ is locally bounded whereas $\textbf{(B3)}$ is true by assuming only that $f_{0,j}$ is locally bounded for each $j\neq 0$.\\
\\
As in Theorem 3.1 in \cite{Carbon-Francq-Tran2010}, the following result gives the asymptotic variance of $f_n$. 
\begin{Prop}\label{variance_de_fn}
Assume that $\textbf{\emph{(A1)}}$ and $\textbf{\emph{(A2)}}$ hold. If one of the following assumptions
\begin{itemize}
\item[$i)$] $\textbf{\emph{(A3)}}$ holds and $\sum_{m\geq 1}m^{2d-1}\,\alpha_{1,1}(m)<\infty$,
\item[$ii)$] $\textbf{\emph{(B3)}}$ holds and $\sum_{m\geq 1}m^{d-1}\,\rho_{1,1}(m)<\infty$,
\end{itemize}
is true then $\lim_{n\to\infty}\vert\Lambda_n\vert b_n\V(f_n(x))=1$ for any $x$ in $\R$ such that $f(x)>0$. 
\end{Prop}
Our main result is the following central limit theorem.
\begin{Th}\label{convergence-loi}
Assume that $\textbf{\emph{(A1)}}$ and $\textbf{\emph{(A2)}}$ hold. If one of the following assumptions
\begin{itemize}
\item[$i)$] $\textbf{\emph{(A3)}}$ holds and $\sum_{m\geq 1}m^{2d-1}\,\alpha_{1,\infty}(m)<\infty$,
\item[$ii)$] $\textbf{\emph{(B3)}}$ holds and $\sum_{m\geq 1}m^{d-1}\,\rho_{1,\infty}(m)<\infty$,
\end{itemize}
is true then for any positive integer $r$ and any distinct points $x_1,...,x_r$ in $\R$ such that $f(x_i)>0$ for any $1\leq i\leq r$,
\begin{equation}\label{limit}
(\vert\Lambda_n\vert b_{n})^{1/2}
\left(\begin{array}{c}
       f_{n}(x_1)-\E f_n(x_1)\\
       \vdots\\
       f_{n}(x_r)-\E f_n(x_r)
       \end{array} \right)
\converge{n}{\infty}{\textrm{$\emph{Law}$}}\mathcal{N}\left(0,\emph{Id}\right)
\end{equation}
where $\emph{Id}$ is the unit matrix of order $r$.
\end{Th}
\textbf{Remark 2}. Theorem \ref{convergence-loi} improves Theorem 4.1 in \cite{Carbon-Francq-Tran2010} in three directions: the regions $\Lambda_n$ where 
the random field is observed are not reduced to rectangular ones (we do not assume any boundary condition), the assumption $\textbf{(A2)}$ on the width bin $b_n$ is minimal and the $\alpha$-mixing condition is weaker than the one assumed in Theorem 4.1 in \cite{Carbon-Francq-Tran2010}, that is 
$\alpha_{1,\infty}(m)=O(m^{-\theta})$ with $\theta>2d$. 
\section{Proofs}
Throughout this section, the symbol $\kappa$ will denote a generic positive constant which the value is not important and 
we recall that $\vert i\vert=\max_{1\leq k\leq d}\vert i_k\vert$ for any $i=(i_1,...,i_d)\in\Z^d$.\\
Let $\tau\in\N^{\ast}\cup\{\infty\}$ be fixed and consider the sequence $(m_{n,\tau})_{n\geq 1}$ defined by 
\begin{equation}\label{def_mn_alpha}
m_{n,\tau}=\max\left\{v_n,\left[\left(\frac{1}{b_n}\sum_{\vert i\vert>v_n}\vert i\vert^d\,\alpha_{1,\tau}(\vert i\vert)\right)^{\frac{1}{d}}\right]+1\right\}
\end{equation}
where $v_n=\big[b_n^{\frac{-1}{2d}}\big]$ and $[\,.\,]$ denotes the integer part function. The following technical lemma is a spatial version of a 
result by Bosq, Merlev\`ede and Peligrad (\cite{Bosq-Merlevede-Peligrad1999}, pages 88-89).
\begin{lemma}\label{mn}
Let $\tau\in\N^{\ast}\cup\{\infty\}$ be fixed. If $\sum_{m\geq 1}m^{2d-1}\,\alpha_{1,\tau}(m)<\infty$ then
$$
m_{n,\tau}\to\infty,\quad m_{n,\tau}^d\,b_n\to0\quad\textrm{and}\quad\frac{1}{m_{n,\tau}^d\,b_n}\sum_{\vert i\vert>m_{n,\tau}}\vert i\vert^d\,\alpha_{1,\tau}(\vert i\vert)\to 0.
$$
\end{lemma}
{\em Proof of Lemma $\ref{mn}$}. First, $m_{n,\tau}$ goes to infinity since $b_n$ goes to zero and $m_{n,\tau}\geq v_n$. We consider the function $\psi$ defined for any $m$ in $\N^{\ast}$ by $\psi(m)=\sum_{\vert i\vert>m}\vert i\vert^d\,\alpha_{1,\tau}(\vert i\vert)$. Since $\sum_{m\geq 1}m^{2d-1}\,\alpha_{1,\tau}(m)<\infty$, we have $\psi(m)$ converges to zero as $m$ goes to infinity. Consequently, $m_{n,\tau}^db_n\leq\max\left\{\sqrt{b_n},\kappa\left(\sqrt{\psi\left(v_n\right)}+b_n\right)\right\}\converge{n}{\infty}{ }0$. Moreover, noting that $m_{n,\tau}^d\geq\frac{1}{b_n}\sqrt{\psi\left(v_n\right)}\geq\frac{1}{b_n}\sqrt{\psi\left(m_{n,\tau}\right)}$ (since $v_n\leq m_{n,\tau}$), we derive also
$$
\frac{1}{m_{n,\tau}^db_n}\sum_{\vert i\vert>m_{n,\tau}}\vert i\vert^d\,\alpha_{1,\tau}(\vert i\vert)\leq\sqrt{\psi(m_{n,\tau})}\converge{n}{\infty}{ }0.
$$
The proof of Lemma $\ref{mn}$ is complete.\\
\\
{\em Proof of Proposition \ref{variance_de_fn}}. For any $n\geq 1$ and any $x$ in $\R$, 
$$
\vert\Lambda_n\vert b_n\V(f_n(x))=\sum_{k\in\Z}\frac{\vert\Lambda_n\vert b_n\V(f_{n,k}(x))}{\sigma^2_{n,k}(x)}\ind{\textrm{J}_{n,k}}(x).
$$
Let $x$ in $\R$ such that $f(x)>0$. For any integer $n\geq 1$, we denote by $k(n)$ the unique integer such that $x$ belongs to $\textrm{J}_{n,k(n)}$. It suffices to show that
$$
\lim_{n\to\infty}\frac{\vert\Lambda_n\vert b_n\V(f_{n,k(n)}(x))}{\sigma^2_{n,k(n)}(x)}=1.
$$
In the sequel, we write $k$ instead of $k(n)$ and we denote $p_s=\int_{\textrm{I}_{n,s}}f(u)du$ for any $s$ in $\Z$. We have 
\begin{align*}
\vert\Lambda_n\vert b_n\V(f_{n,k}(x))
&=a^2_k(x)\left(\frac{p_k(1-p_k)}{b_n}+\frac{1}{\vert\Lambda_n\vert b_n}\sum_{\stackrel{i,j\in\Lambda_n}{i\neq j}}\textrm{Cov}(Y_{i,k},Y_{j,k})\right)\\
&+\overline{a}^2_k(x)\left(\frac{p_{k+1}(1-p_{k+1})}{b_n}+\frac{1}{\vert\Lambda_n\vert b_n}\sum_{\stackrel{i,j\in\Lambda_n}{i\neq j}}\textrm{Cov}(Y_{i,k+1},Y_{j,k+1})\right)\\
&+2a_k(x)\overline{a}_k(x)\left(\frac{-p_{k}p_{k+1}}{b_n}+\frac{1}{\vert\Lambda_n\vert b_n}\sum_{\stackrel{i,j\in\Lambda_n}{i\neq j}}\textrm{Cov}(Y_{i,k},Y_{j,k+1})\right).
\end{align*}
Denote also
$$
w_{n,k}(x)=\frac{1}{b_n}\left(a^2_k(x)p_k(1-p_k)+\overline{a}^2_k(x)p_{k+1}(1-p_{k+1})-2a_k(x)\overline{a}_k(x)p_{k}p_{k+1}\right).
$$
Arguing as in the proof of Lemma 3.2 in \cite{Carbon-Francq-Tran2010}, we have by Taylor expansion
\begin{align*}
p_k&=b_nf(x)+\frac{b_n}{2}(2(kb_n-x)-b_n)f^{'}(c_k)\\
p_{k+1}&=b_nf(x)+\frac{b_n}{2}(2(kb_n-x)+b_n)f^{'}(c_{k+1})
\end{align*}
where $c_k\in \textrm{J}_{n,k}$ and $c_{k+1}\in \textrm{J}_{n,k+1}$. Then for $j=k$ or $j=k+1$,
\begin{equation}\label{encadrement_p_j}
\max\{0,f(x)b_n-\kappa b_n^2\}\leq p_{j}\leq f(x)b_n+\kappa b_n^2.
\end{equation}
Consequently
\begin{equation}\label{p_k_p_k+1}
\max\{0,b_nf^2(x)-2b_n^2\kappa f(x)+\kappa^2b_n^3\}\leq\frac{p_kp_{k+1}}{b_n}\leq b_nf^2(x)+2b_n^2\kappa f(x)-\kappa^2b_n^3
\end{equation}
and for $j=k$ or $j=k+1$,
\begin{equation}\label{p_j_1_moins_p_j}
\max\{0,f(x)-(\kappa+f^2(x)b_n)+\kappa^2b_n^3\}\leq\frac{p_j(1-p_j)}{b_n}\leq f(x)+(\kappa-f^2(x))b_n+\kappa^2b_n^3.
\end{equation}
Finally, we obtain
$$
\lim_{n\to\infty}\frac{w_{n,k}(x)}{\sigma^2_{n,k}(x)}=1.
$$
Let $(s,t)$ be equal to $(k,k)$, $(k,k+1)$ or $(k+1,k+1)$. It suffices to show that 
\begin{equation}\label{limit_covariance}
\lim_{n\to\infty}\frac{1}{\vert\Lambda_n\vert b_n}\sum_{\stackrel{i,j\in\Lambda_n}{i\neq j}}\textrm{Cov}(Y_{i,s},Y_{j,t})=0
\end{equation}
By stationarity of the random field $(X_i)_{i\in\Z^d}$, we have
\begin{equation}\label{inequality_sum_covariance_Y}
\frac{1}{\vert\Lambda_n\vert b_n}\left\vert\sum_{\stackrel{i,j\in\Lambda_n}{i\neq j}}\textrm{Cov}(Y_{i,s},Y_{j,t})\right\vert
\leq\frac{1}{b_n}\sum_{\stackrel{j\in\Z^d}{j\neq 0}}\vert\textrm{Cov}(Y_{0,s},Y_{j,t})\vert.
\end{equation}
Using (\ref{encadrement_p_j}), for $j\neq 0$, we have
\begin{equation}\label{rho_mixing_inequality_0}
\vert\textrm{Cov}(Y_{0,s},Y_{j,t})\vert\leq \|Y_{0,s}\|_2\|Y_{0,t}\|_2\rho_{1,1}(\vert j\vert)=\sqrt{p_sp_t}\rho_{1,1}(\vert j\vert)\leq\kappa b_n\rho_{1,1}(\vert j\vert).
\end{equation}
Moreover, using again (\ref{encadrement_p_j}), for any $j\neq 0$, 
\begin{equation}\label{covariance_Y}
\vert\textrm{Cov}(Y_{0,s},Y_{j,t})\vert\leq\P\left(X_0\in I_{n,s}, X_j\in I_{n,t}\right)+\kappa b_n^2.
\end{equation}
Assuming $\textbf{(B3)}$ and $\sum_{m\geq 1}m^{d-1}\rho_{1,1}(m)<\infty$ and combining (\ref{rho_mixing_inequality_0}) and (\ref{covariance_Y}) with the dominated convergence theorem, we obtain 
\begin{equation}\label{CDL}
\sum_{\stackrel{j\in\Z^d}{j\neq 0}}\frac{\vert\textrm{Cov}(Y_{0,s},Y_{j,t})\vert}{b_n}=o(1).
\end{equation}
Finally, (\ref{limit_covariance}) follows from inequality (\ref{inequality_sum_covariance_Y}). Similarly, by Rio's covariance inequality (cf. \cite{Rio1993}, Theorem 1.1),
$$
\vert\textrm{Cov}(Y_{0,s},Y_{j,t})\vert\leq 2\int_0^{2\alpha\left(\sigma(Y_{0,s}),\sigma(Y_{j,t})\right)}Q_{Y_{0,s}}(u)Q_{Y_{j,t}}(u)du
$$
where $Q_{Z}(u)=\inf\{t;\,\P(\vert Z\vert>t)\leq u\}$ for any $u$ in $[0,1]$. Since $Y_{0,s}$ and $Y_{j,t}$ are bounded by $1$, we derive
\begin{equation}\label{covariance_Ybis}
\vert\textrm{Cov}(Y_{0,s},Y_{j,t})\vert\leq 4\alpha_{1,1}(\vert j\vert).
\end{equation}
Using (\ref{covariance_Y}) and assuming $\textbf{(A3)}$, we derive
\begin{equation}\label{covariance_Y_alpha}
\sup_{j\neq 0}\vert\textrm{Cov}(Y_{0,s},Y_{j,t})\vert\leq\kappa b_n^2.
\end{equation}
Assuming $\sum_{m\geq 1}m^{2d-1}\alpha_{1,1}(m)<\infty$ and combining (\ref{covariance_Ybis}), (\ref{covariance_Y_alpha}) and Lemma \ref{mn}, we obtain
\begin{equation}\label{CDL_alpha}
\sum_{\stackrel{j\in\Z^d}{j\neq 0}}\frac{\vert\textrm{Cov}(Y_{0,s},Y_{j,t})\vert}{b_n}\leq \kappa\left(m_{n,1}^db_n+\frac{1}{m_{n,1}^db_n}\sum_{\vert j\vert>m_{n,1}}\vert j\vert^d\alpha_{1,1}(\vert j\vert)\right)=o(1)
\end{equation}
where $(m_{n,1})_{n\geq 1}$ is defined by (\ref{def_mn_alpha}). Finally, using inequality (\ref{inequality_sum_covariance_Y}), we obtain again (\ref{limit_covariance}). The proof of Proposition \ref{variance_de_fn} is complete.\\
\\
\textbf{Remark 3}. The reader should note that the asymptotic variance given in Theorem 3.1 in \cite{Carbon-Francq-Tran2010} is not the good one. 
In fact, using the notations in \cite{Carbon-Francq-Tran2010}, it should be $(1/2+2(k_0-x/b)^2)f(x)$ instead of $(1/2+(2k_0-x/b)^2)f(x)$.\\
\\
{\em Proof of Theorem \ref{convergence-loi}}. Without loss of generality, we assume that $r=2$ and we denote $x$ and $y$ in place of $x_1$ and $x_2$.
Let $\lambda_1$ and $\lambda_2$ be fixed in $\R$ such that $\lambda_1^2+\lambda_2^2=1$. For any $i$ in $\Z^d$, we define 
$\Delta_i=\lambda_1Z_{i}(x)+\lambda_2 Z_{i}(y)$ where for any $u$ in $\R$ such that $f(u)>0$,
$$
Z_i(u)=\frac{1}{\sqrt{b_n}}\sum_{s\in\Z}\frac{a_s(u)(Y_{i,s}-\E(Y_{i,s}))+\overline{a}_s(u)(Y_{i,s+1}-\E(Y_{i,s+1}))}{\sigma_{n,s}(u)}\ind{\textrm{J}_{n,s}}(u),
$$
$a_s(u)=\frac{1}{2}+s-\frac{u}{b_n},\,\,\,\overline{a}_s(u)=a_{-s}(-u)\,\,\,\textrm{and}\,\,\,\sigma^2_{n,s}(u)=\left(\frac{1}{2}+2\left(s-\frac{u}{b_n}\right)^2\right)f(u)$.
\begin{lemma}\label{esperance_delta0_deltai}
$\E(\Delta_0^2)\converge{n}{\infty}{ }1$ and $\sum_{i\in W}\E\vert\Delta_0\Delta_i\vert\leq\kappa\vert W\vert b_n+o(1)$ for any finite subset $W$ of $\Z^d\backslash\{0\}$.  
\end{lemma}
{\em Proof of Lemma \ref{esperance_delta0_deltai}}. Recall that $p_s=\int_{\textrm{I}_{n,s}}f(u)du$ 
and $\textrm{J}_{n,s}=[(s-\frac{1}{2})b_n,(s+\frac{1}{2})b_n)$ for any $s$ in $\Z$. We have 
\begin{equation}\label{Delta_0_carre}
\E(\Delta_0^2)=\lambda_1^2\E(Z_{0}^2(x))+\lambda_2^2\E(Z_{0}^2(y))+2\lambda_1\lambda_2\E(Z_{0}(x)Z_{0}(y))
\end{equation}
and 
\begin{align*}
b_n\E(Z_{0}(x)Z_{0}(y))
&=\sum_{(k,l)\in\Z^2}\frac{a_k(x)a_l(y)}{\sigma_{n,k}(x)\sigma_{n,l}(y)}\E(Y_{0,k}-p_k)(Y_{0,l}-p_l)\ind{\textrm{J}_{n,k}\times \textrm{J}_{n,l}}(x,y)\\
&+\sum_{(k,l)\in\Z^2}\frac{a_k(x)\overline{a}_l(y)}{\sigma_{n,k}(x)\sigma_{n,l}(y)}\E(Y_{0,k}-p_k)(Y_{0,l+1}-p_{l+1})\ind{\textrm{J}_{n,k}\times \textrm{J}_{n,l}}(x,y)\\
&+\sum_{(k,l)\in\Z^2}\frac{\overline{a}_k(x)a_l(y)}{\sigma_{n,k}(x)\sigma_{n,l}(y)}\E(Y_{0,k+1}-p_{k+1})(Y_{0,l}-p_l)\ind{\textrm{J}_{n,k}\times \textrm{J}_{n,l}}(x,y)\\
&+\sum_{(k,l)\in\Z^2}\frac{\overline{a}_k(x)\overline{a}_l(y)}{\sigma_{n,k}(x)\sigma_{n,l}(y)}\E(Y_{0,k+1}-p_{k+1})(Y_{0,l+1}-p_{l+1})\ind{\textrm{J}_{n,k}\times \textrm{J}_{n,l}}(x,y)
\end{align*}
For $n$ sufficiently large, $(x,y)$ belongs to $\textrm{J}_{n,k}\times \textrm{J}_{n,l}$ with $\vert k-l\vert\geq 2$. Then for any $s=k$ or $s=k+1$ and $t=l$ or $t=l+1$, 
$\vert\E(Y_{0,s}-p_s)(Y_{0,t}-p_t)\vert=p_sp_t\leq \kappa b_n^2$. Since $0\leq a_s(u)\leq 1$, $0\leq \overline{a}_s(u)\leq 1$ and $\sigma^2_{n,s}(u)\geq f(u)/2>0$ 
for any $u$ in $\textrm{J}_{n,s}$ such that $f(u)>0$ and any $s$ in $\Z$, we obtain 
\begin{equation}\label{limite_esperance_Z_0_x_Z_0_y}
\left\vert\E(Z_{0}(x)Z_{0}(y))\right\vert\leq \kappa b_n.
\end{equation}
Similarly, for any $u$ in $\R$, we have
\begin{align*}
b_n\E(Z^2_{0}(u))
&=\sum_{k\in\Z}\frac{a^2_k(u)}{\sigma^2_{n,k}(u)}\E(Y_{0,k}-p_k)^2\ind{\textrm{J}_{n,k}}(u)\\
&+2\sum_{k\in\Z}\frac{a_k(u)\overline{a}_k(u)}{\sigma^2_{n,k}(u)}\E(Y_{0,k}-p_k)(Y_{0,k+1}-p_{k+1})\ind{\textrm{J}_{n,k}}(u)\\
&+\sum_{k\in\Z}\frac{\overline{a}^2_k(u)}{\sigma^2_{n,k}(u)}\E(Y_{0,k+1}-p_{k+1})^2\ind{\textrm{J}_{n,k}}(u).
\end{align*}
Noting that $\E(Y_{0,s}-p_s)^2=p_s(1-p_s)$ and $\E(Y_{0,s}-p_s)(Y_{0,s+1}-p_{s+1})=-p_sp_{s+1}$ for any $s$ in $\Z$ and 
keeping in mind ($\ref{p_k_p_k+1}$) and ($\ref{p_j_1_moins_p_j}$), we obtain for any $u$ in $\R$,
\begin{equation}\label{esperance_Z_0_carre}
\E(Z^2_{0}(u))\converge{n}{\infty}{}1.
\end{equation}
Combining ($\ref{Delta_0_carre}$), ($\ref{limite_esperance_Z_0_x_Z_0_y}$) and ($\ref{esperance_Z_0_carre}$), we obtain 
$\E(\Delta_0^2)\converge{n}{\infty}{}\lambda_1^2+\lambda_2^2=1$.\\
\\
Let $W$ be a finite subset of $\Z^d\backslash\{0\}$ and let $i\in W$ be fixed. We have
\begin{equation}\label{majoration_delta_0_delta_i}
\E\vert\Delta_0\Delta_i\vert\leq\lambda_1^2\E\vert Z_{0}(x)Z_{i}(x)\vert+2\vert\lambda_1\vert\vert\lambda_2\vert\E\vert Z_{0}(x)Z_{i}(y)\vert
+\lambda_2^2\E\vert Z_{0}(y)Z_{i}(y)\vert. 
\end{equation}
If $u$ and $v$ are fixed in $\R$ then
\begin{align*}
b_n\E\vert Z_{0}(u)Z_{i}(v)\vert
&\leq\sum_{(k,l)\in\Z^2}\frac{\vert\overline{a}_k(u)\overline{a}_{l}(v)\vert}{\sigma_{n,k}(u)\sigma_{n,l}(v)}\E\big\vert(Y_{0,k}-p_k)(Y_{i,l}-p_l)\big\vert\ind{\textrm{J}_{n,k}\times \textrm{J}_{n,l}}(u,v)\\
&\quad+\sum_{(k,l)\in\Z^2}\frac{\vert\overline{a}_k(u)a_l(v)\vert}{\sigma_{n,k}(u)\sigma_{n,l}(v)}\E\big\vert(Y_{0,k}-p_k)(Y_{i,l+1}-p_{l+1})\big\vert\ind{\textrm{J}_{n,k}\times \textrm{J}_{n,l}}(u,v)\\
&\quad+\sum_{(k,l)\in\Z^2}\frac{\vert\overline{a}_k(u)\overline{a}_l(v)\vert}{\sigma_{n,k}(u)\sigma_{n,l}(v)}\E\big\vert(Y_{0,k+1}-p_{k+1})(Y_{i,l}-p_l)\big\vert\ind{\textrm{J}_{n,k}\times \textrm{J}_{n,l}}(u,v)\\
&\quad+\sum_{(k,l)\in\Z^2}\frac{\vert a_k(u)a_l(v)\vert}{\sigma_{n,k}(u)\sigma_{n,l}(v)}\E\big\vert(Y_{0,k+1}-p_{k+1})(Y_{i,l+1}-p_{l+1})\big\vert\ind{\textrm{J}_{n,k}\times \textrm{J}_{n,l}}(u,v).
\end{align*}
Noting that for any $s$ and $t$ in $\Z$,
$$
\E\vert(Y_{0,s}-p_s)(Y_{i,t}-p_t)\vert\leq\E(Y_{0,s} Y_{i,t})+3p_sp_t\leq \vert\textrm{Cov}(Y_{0,s},Y_{i,t})\vert+\kappa b_n^2,
$$
we derive
\begin{align*}
\E\vert Z_{0}(u)Z_{i}(v)\vert
&\leq\sum_{(k,l)\in\Z^2}\frac{\vert\overline{a}_k(u)\overline{a}_{l}(v)\vert}{\sigma_{n,k}(u)\sigma_{n,l}(v)}\frac{\vert\textrm{Cov}(Y_{0,k},Y_{i,l})\vert}{b_n}\ind{\textrm{J}_{n,k}\times \textrm{J}_{n,l}}(u,v)\\
&\quad+\sum_{(k,l)\in\Z^2}\frac{\vert\overline{a}_k(u)a_l(v)\vert}{\sigma_{n,k}(u)\sigma_{n,l}(v)}\frac{\vert\textrm{Cov}(Y_{0,k},Y_{i,l+1})\vert}{b_n}\ind{\textrm{J}_{n,k}\times \textrm{J}_{n,l}}(u,v)\\
&\quad+\sum_{(k,l)\in\Z^2}\frac{\vert\overline{a}_k(u)\overline{a}_l(v)\vert}{\sigma_{n,k}(u)\sigma_{n,l}(v)}\frac{\vert\textrm{Cov}(Y_{0,k+1},Y_{i,l})\vert}{b_n}\ind{\textrm{J}_{n,k}\times \textrm{J}_{n,l}}(u,v)\\
&\quad+\sum_{(k,l)\in\Z^2}\frac{\vert a_k(u)a_l(v)\vert}{\sigma_{n,k}(u)\sigma_{n,l}(v)}\frac{\vert\textrm{Cov}(Y_{0,k+1},Y_{i,l+1})\vert}{b_n}\ind{\textrm{J}_{n,k}\times \textrm{J}_{n,l}}(u,v)\\
&\quad+\kappa b_n.
\end{align*}
Assuming $\textbf{(B3)}$ and $\sum_{m>0}m^{d-1}\rho_{1,1}(m)<\infty$ and using (\ref{CDL}) or assuming $\textbf{(A3)}$ and $\sum_{m>0}m^{2d-1}\alpha_{1,1}(m)<\infty$ and using (\ref{CDL_alpha}), we obtain 
\begin{equation}\label{borne_sum_E_Z0Zi}
\sum_{i\in W}\E\vert Z_{0}(u)Z_{i}(v)\vert \leq \kappa\vert W\vert b_n+o(1).
\end{equation}
Combining $(\ref{majoration_delta_0_delta_i})$ and ($\ref{borne_sum_E_Z0Zi}$), we obtain $\sum_{i\in W}\E\vert \Delta_{0}\Delta_{i}\vert \leq \kappa\vert W\vert b_n+o(1)$.
The proof of Lemma \ref{esperance_delta0_deltai} is complete.\\
\\
We are going to follow the Lindeberg-type proof of the central limit theorem for stationary random fields established in \cite{Dedecker1998}. 
Let $\varphi$ be a one to one map from $[1,\kappa]\cap\N^{\ast}$ to a finite subset of $\Z^d$ and $(\xi_i)_{i\in\Z^d}$ a real random field. For
all integers $k$ in $[1,\kappa]$, we denote
$$
S_{\varphi(k)}(\xi)=\sum_{i=1}^k \xi_{\varphi(i)}\quad\textrm{and}\quad
S_{\varphi(k)}^{c}(\xi)=\sum_{i=k}^\kappa \xi_{\varphi(i)}
$$
with the convention $S_{\varphi(0)}(\xi)=S_{\varphi(\kappa+1)}^{c}(\xi)=0$. On the lattice $\Z^{d}$ we define the lexicographic order as follows: if
$i=(i_{1},...,i_{d})$ and $j=(j_{1},...,j_{d})$ are distinct elements of $\Z^{d}$, the notation $i<_{\textrm{lex}}j$ means that either
$i_{1}<j_{1}$ or for some $p$ in $\{2,3,...,d\}$, $i_{p}<j_{p}$ and $i_{q}=j_{q}$ for $1\leq q<p$. To
describe the set $\Lambda_{n}$, we define the one to one map $\varphi$ from $[1,\vert\Lambda_n\vert]\cap\N^{\ast}$ to $\Lambda_{n}$ by:
$\varphi$ is the unique function such that $\varphi(k)<_{\text{\text{lex}}}\varphi(l)$ for $1\leq k<l\leq\vert\Lambda_n\vert$. 
From now on, we consider a field $(\xi_{i})_{i\in\Z^d}$ of i.i.d. random variables independent of
$(X_{i})_{i\in\Z^d}$ such that $\xi_{0}$ has the standard normal
law $\mathcal{N}(0,1)$. We introduce the fields $\Gamma$ and $\gamma$ defined for any $i$ in $\Z^d$ by 
$$
\Gamma_{i}=\frac{\Delta_i}{\vert\Lambda_n\vert^{1/2}}\quad\textrm{and}\quad\gamma_{i}=\frac{\xi_i}{\vert\Lambda_n\vert^{1/2}}
$$ 
Let $h$ be any function from $\R$ to $\R$. For $0\leq k\leq l\leq\vert\Lambda_n\vert+1$, we introduce
$h_{k,l}(\Gamma)=h(S_{\varphi(k)}(\Gamma)+S_{\varphi(l)}^{c}(\gamma))$. With the above
convention we have that $h_{k,\vert\Lambda_n\vert+1}(\Gamma)=h(S_{\varphi(k)}(\Gamma))$ and also
$h_{0,l}(\Gamma)=h(S_{\varphi(l)}^{c}(\gamma))$. In the sequel, we will often
write $h_{k,l}$ instead of $h_{k,l}(\Gamma)$. We denote by $B_{1}^4(\R)$ the unit ball of $C_{b}^4(\R)$: $h$ belongs to $B_{1}^4(\R)$ if and only if 
$h$ is bounded by $1$, belongs to $C^4(\R)$ and its first four derivatives are also bounded by $1$. It suffices to prove that for all $h$ in $B_{1}^4(\R)$,
$$
\E\left(h\left(S_{\varphi(\vert\Lambda_n\vert)}(\Gamma)\right)\right)\converge{n}{\infty}{}\E \left(h\left(\xi_0\right)\right).
$$
We use Lindeberg's decomposition:
$$
\E\left(h\left(S_{\varphi(\vert\Lambda_n\vert)}(\Gamma)\right)-h\left(\xi_{0}\right)\right)
=\sum_{k=1}^{\vert\Lambda_n\vert}\E \left(h_{k,k+1}-h_{k-1,k}\right).
$$
Now,
$$
h_{k,k+1}-h_{k-1,k}=h_{k,k+1}-h_{k-1,k+1}+h_{k-1,k+1}-h_{k-1,k}.
$$
Applying Taylor's formula we get that:
$$
h_{k,k+1}-h_{k-1,k+1}=\Gamma_{\varphi(k)}h_{k-1,k+1}^{'}+\frac{1}{2}\Gamma_{\varphi(k)}^{2}h_{k-1,k+1}^{''}+R_{k}
$$
and
$$
h_{k-1,k+1}-h_{k-1,k}=-\gamma_{\varphi(k)}h_{k-1,k+1}^{'}-\frac{1}{2}\gamma_{\varphi(k)}^{2}h_{k-1,k+1}^{''}+r_{k}
$$
where $\vert R_{k}\vert\leq \Gamma_{\varphi(k)}^2(1\wedge\vert \Gamma_{\varphi(k)}\vert)$ and $\vert r_{k}\vert\leq\gamma_{\varphi(k)}^2(1\wedge\vert\gamma_{\varphi(k)}\vert)$.
Since $(\Gamma,\xi_{i})_{i\neq \varphi(k)}$ is independent of $\xi_{\varphi(k)}$, it follows that
$$
\E \left(\gamma_{\varphi(k)}h_{k-1,k+1}^{'}\right)=0\quad\textrm{and}\quad
\E \left(\gamma_{\varphi(k)}^2h_{k-1,k+1}^{''}\right)=\E \left(\frac{h_{k-1,k+1}^{''}}{\vert\Lambda_n\vert}\right)
$$
Hence, we obtain
\begin{align*}
\E \left(h(S_{\varphi(\vert\Lambda_n\vert)}(\Gamma))-h\left(\xi_0\right)\right)&=
\sum_{k=1}^{\vert\Lambda_n\vert}\E (\Gamma_{\varphi(k)}h_{k-1,k+1}^{'})\\
&\quad+\sum_{k=1}^{\vert\Lambda_n\vert}\E \left(\left(\Gamma_{\varphi(k)}^2-\frac{1}{\vert\Lambda_n\vert}\right)\frac{h_{k-1,k+1}^{''}}{2}\right)\\
&\quad+\sum_{k=1}^{\vert\Lambda_n\vert}\E \left(R_{k}+r_{k}\right).
\end{align*}
Let $1\leq k\leq \vert\Lambda_n\vert$ be fixed. Since $\vert\Delta_0\vert$ is bounded by $\kappa/\sqrt{b_n}$, applying Lemma $\ref{esperance_delta0_deltai}$
we derive
$$
\E\vert R_k\vert\leq\frac{\E\vert\Delta_0\vert^3}{\vert\Lambda_n\vert^{3/2}}\leq\frac{\kappa}{(\vert\Lambda_n\vert^{3}\,b_n)^{1/2}}
\quad\textrm{and}\quad
\E\vert r_k\vert\leq\frac{\E\vert\xi_0\vert^3}{\vert\Lambda_n\vert^{3/2}}.
$$
Consequently, we obtain
$$
\sum_{k=1}^{\vert\Lambda_n\vert}\E \left(\vert R_{k}\vert+\vert r_{k}\vert\right)=O\left(\frac{1}{(\vert\Lambda_n\vert b_n)^{1/2}}+\frac{1}{\vert\Lambda_n\vert^{1/2}}\right)=o(1).
$$
Now, it is sufficient to show
\begin{equation}\label{equation1}
\lim_{n\to\infty}\sum_{k=1}^{\vert\Lambda_n\vert}\left(\E (\Gamma_{\varphi(k)}h_{k-1,k+1}^{'})+\E \left(\left(\Gamma_{\varphi(k)}^2-\frac{1}{\vert\Lambda_n\vert}\right)\frac{h_{k-1,k+1}^{''}}{2}\right)\right)=0.
\end{equation}
First, we focus on $\sum_{k=1}^{\vert\Lambda_n\vert}\E \left(\Gamma_{\varphi(k)}h_{k-1,k+1}^{'}\right)$. Let the sets 
$\{V_{i}^{k}\,;\,i\in\Z^{d}\,,\,k\in\N^{\ast}\}$ be defined as follows: $V_{i}^{1}=\{j\in\Z^{d}\,;\,j<_{\textrm{lex}}i\}$ and 
$V_{i}^{k}=V_{i}^{1}\cap\{j\in\Z^{d}\,;\,\vert i-j\vert\geq k\}$ for $k\geq 2$. Let $(N_n)_{n\geq 1}$ be a sequence of positive integers satisfying
\begin{equation}\label{Nn}
N_n\to\infty\quad\textrm{such that}\quad N_n^db_n\to0.
\end{equation}
For all $n$ in $\N^{\ast}$ and all integer $1\leq k\leq\vert\Lambda_n\vert$, we define
$$
\textrm{E}_{k}^{(n)}=\varphi([1,k]\cap\N^{\ast})\cap V_{\varphi(k)}^{N_n}\quad\textrm{and}\quad
S_{\varphi(k)}^{(n)}(\Gamma)=\sum_{i\in\textrm{E}_{k}^{(n)}}\Gamma_{i}.
$$ 
For any function $\Psi$ from $\R$ to $\R$, we define $\Psi_{k-1,l}^{(n)}=\Psi(S_{\varphi(k)}^{(n)}(\Gamma)+S_{\varphi(l)}^c(\gamma))$. We are going to 
apply this notation to the successive derivatives of the function $h$. Our aim is to show that
\begin{equation}\label{equation1bis}
\lim_{n\to\infty}\sum_{k=1}^{\vert\Lambda_n\vert}\E\left(\Gamma_{\varphi(k)}h_{k-1,k+1}^{'}-\Gamma_{\varphi(k)}\left(S_{\varphi(k-1)}(\Gamma)-S_{\varphi(k)}^{(n)}(\Gamma)\right)h_{k-1,k+1}^{''}\right)=0.
\end{equation}
First, we use the decomposition
$$
\Gamma_{\varphi(k)}h_{k-1,k+1}^{'}=\Gamma_{\varphi(k)}h_{k-1,k+1}^{'(n)}+\Gamma_{\varphi(k)}\left(h_{k-1,k+1}^{'}-h_{k-1,k+1}^{'(n)}\right).
$$
We consider a one to one map $\psi$ from $[1,\vert \textrm{E}_{k}^{(n)}\vert]\cap\N^{\ast}$ to $\textrm{E}_{k}^{(n)}$ such that $\vert
\psi(i)-\varphi(k)\vert\leq\vert \psi(i-1)-\varphi(k)\vert$. For any subset $B$ of $\Z^{d}$, recall that $\F_{B}=\sigma(X_{i}\,;\,i\in B)$ and set
$$
\E_{M}(X_{i})=\E(X_{i}\vert\F_{V_{i}^{M}}),\quad M\in \N^{\ast},\quad i\in\Z^d.
$$
The choice of the map $\psi$ ensures that $S_{\psi(i)}(\Gamma)$ and $S_{\psi(i-1)}(\Gamma)$ are $\F_{V_{\varphi(k)}^{\vert \psi(i)-\varphi(k)\vert}}$-measurable.
The fact that $\gamma$ is independent of $\Gamma$ imply that $\E \left(\Gamma_{\varphi(k)}h^{'}\left(S_{\varphi(k+1)}^c(\gamma)\right)\right)=0$. Therefore
\begin{equation}\label{equation_theta}
\left\vert \E \left(\Gamma_{\varphi(k)}h_{k-1,k+1}^{'{(n)}}\right)\right\vert=\left\vert\sum_{i=1}^{\left\vert
E_{k}^{(n)}\right\vert}\E \left(\Gamma_{\varphi(k)}\left(\theta_i-\theta_{i-1}\right)\right)\right\vert
\end{equation}
where $\theta_i=h^{'}\left(S_{\psi(i)}(\Gamma)+S_{\varphi(k+1)}^c(\gamma)\right)$. Since $S_{\psi(i)}(\Gamma)$ and $S_{\psi(i-1)}(\Gamma)$ are $\F_{V_{\varphi(k)}^{\vert
\psi(i)-\varphi(k)\vert}}$-measurable, we can take the conditional expectation of $\Gamma_{\varphi(k)}$ with respect to $\F_{V_{\varphi(k)}^{\vert \psi(i)-\varphi(k)\vert}}$ 
in the right hand side of ($\ref{equation_theta}$). On the other hand the function $h^{'}$ is $1$-Lipschitz, hence 
$\left\vert\theta_i-\theta_{i-1}\right\vert\leq\vert \Gamma_{\psi(i)}\vert$. Consequently, 
$\left\vert\E \left(\Gamma_{\varphi(k)}\left(\theta_i-\theta_{i-1}\right)\right)\right\vert\leq\E \vert \Gamma_{\psi(i)}\E _{\vert \psi(i)-\varphi(k)\vert}\left(\Gamma_{\varphi(k)}\right)\vert$ 
and
$$
\left\vert
\E \left(\Gamma_{\varphi(k)}h_{k-1,k+1}^{'(n)}\right)\right\vert\leq\sum_{i=1}^{\left\vert E_{k}^{(n)}\right\vert}
\E \vert \Gamma_{\psi(i)}\E _{\vert \psi(i)-\varphi(k)\vert}(\Gamma_{\varphi(k)})\vert.
$$
Hence,
\begin{align*}
\left\vert\sum_{k=1}^{\vert\Lambda_n\vert}\E \left(\Gamma_{\varphi(k)}h_{k-1,k+1}^{'(n)}\right)\right\vert
&\leq\frac{1}{\vert\Lambda_n\vert}\sum_{k=1}^{\vert\Lambda_n\vert}\sum_{i=1}^{\vert E_{k}^{(n)}\vert} \E \vert\Delta_{\psi(i)}\E _{\vert \psi(i)-\varphi(k)\vert}(\Delta_{\varphi(k)})\vert\\
&\leq \sum_{\vert j\vert \geq N_n}\|\Delta_j\E _{\vert j\vert}(\Delta_0)\|_1.
\end{align*}
For any $j$ in $\Z^d$, we have $\|\Delta_j\E _{\vert j\vert}(\Delta_0)\|_1
=\textrm{Cov}\left(\vert\Delta_j\vert\left(\ind{\E_{\vert j\vert}(\Delta_0)\geq 0}-\ind{\E_{\vert j\vert}(\Delta_0)<0}\right),\Delta_0\right)$ and consequently $\|\Delta_j\E _{\vert j\vert}(\Delta_0)\|_1\leq \|\Delta_0\|_2^2\rho_{1,\infty}(\vert j\vert)$. By Lemma \ref{esperance_delta0_deltai}, we know that $\|\Delta_0\|_2$ is bounded. So, assuming $\sum_{m\geq 0}m^{d-1}\rho_{1,\infty}(m)<\infty$, we derive
\begin{equation}\label{petit_o_de_1}
\left\vert\sum_{k=1}^{\vert\Lambda_n\vert}\E \left(\Gamma_{\varphi(k)}h_{k-1,k+1}^{'(n)}\right)\right\vert
=o(1).
\end{equation}
By Rio's covariance inequality (cf. \cite{Rio1993}, Theorem 1.1), we have also 
$$
\|\Delta_j\E _{\vert j\vert}(\Delta_0)\|_1\leq 4\int_{0}^{\alpha_{1,\infty}(\vert j\vert)}Q_{\Delta_0}^2(u)du
$$
where $Q_{\Delta_0}$ is defined by $Q_{\Delta_0}(u)=\inf\{t\geq 0\,;\,\P(\vert \Delta_0\vert>t)\leq u\}$ for any $u$ in $[0,1]$. Since $\vert\Delta_0\vert$ is bounded by 
$\kappa/\sqrt{b_n}$, we have $Q_{\Delta_0}\leq \kappa/\sqrt{b_n}$ and $\|\Delta_j\E _{\vert j\vert}(\Delta_0)\|_1\leq\frac{\kappa}{b_n}\,\alpha_{1,\infty}(\vert j\vert)$. 
Assuming $\sum_{m\geq 0}m^{2d-1}\alpha_{1,\infty}(m)<\infty$ and choosing $N_n=m_{n,\infty}$ (recall that $(m_{n,\infty})_{n\geq 1}$ is defined by (\ref{def_mn_alpha}) and satisfies $m_{n,\infty}\to\infty$ such that $m_{n,\infty}^db_n\to 0$), we derive
$$
\left\vert\sum_{k=1}^{\vert\Lambda_n\vert}\E \left(\Gamma_{\varphi(k)}h_{k-1,k+1}^{'(n)}\right)\right\vert
\leq\frac{\kappa}{m_{n,\infty}^db_n}\sum_{\vert j\vert \geq m_{n,\infty}}\vert j\vert^d\,\alpha_{1,\infty}(\vert j\vert).
$$
By Lemma \ref{mn}, we obtain again (\ref{petit_o_de_1}). Now, applying Taylor's formula,
$$
\Gamma_{\varphi(k)}(h_{k-1,k+1}^{'}-h_{k-1,k+1}^{'(n)})=\Gamma_{\varphi(k)}(S_{\varphi(k-1)}(\Gamma)-S_{\varphi(k)}^{(n)}(\Gamma))h_{k-1,k+1}^{''}+R_{k}^{'},
$$
where $\vert R_{k}^{'}\vert\leq 2\vert \Gamma_{\varphi(k)}(S_{\varphi(k-1)}(\Gamma)-S_{\varphi(k)}^{(n)}(\Gamma))(1\wedge\vert S_{\varphi(k-1)}(\Gamma)-S_{\varphi(k)}^{(n)}(\Gamma)\vert)\vert$. 
Consequently $(\ref{equation1bis})$ holds if and only if $\lim_{n\to\infty}\sum_{k=1}^{\vert\Lambda_n\vert}\E \vert R_{k}^{'}\vert=0$. In fact, 
denoting $W_n=\{-N_n+1,...,N_n-1\}^d$ and $W_n^{\ast}=W_n\backslash\{0\}$, we have
\begin{align*}
\sum_{k=1}^{\vert\Lambda_n\vert}\E \vert R_{k}^{'}\vert 
&\leq 2\E \left(\vert\Delta_{0}\vert\left(\sum_{i\in W_n}\vert\Delta_{i}\vert\right)
\left(1\wedge\frac{1}{\vert\Lambda_n\vert^{1/2}}\sum_{i\in W_n}\vert\Delta_{i}\vert\right)\right)\\
&=2\E\left(\left(\Delta_0^2+\sum_{i\in W_n^{\ast}}\vert\Delta_0\Delta_i\vert\right)\left(1\wedge\frac{1}{\vert\Lambda_n\vert^{1/2}}\sum_{i\in W_n}\vert\Delta_i\vert\right)\right)\\
&\leq\frac{2}{\vert\Lambda_n\vert^{1/2}}\sum_{i\in W_n}\E(\Delta_0^2\vert\Delta_i\vert)+2\sum_{i\in W_n^{\ast}}\E\vert\Delta_0\Delta_i\vert.
\end{align*}
Since $\vert\Delta_0\vert\leq\frac{\kappa}{\sqrt{b_n}}$, we derive
\begin{align*}
\sum_{k=1}^{\vert\Lambda_n\vert}\E \vert R_{k}^{'}\vert
&\leq\frac{\kappa}{(\vert\Lambda_n\vert b_n)^{1/2}}\sum_{i\in W_n}\E(\vert\Delta_0\Delta_i\vert)+2\sum_{i\in W_n^{\ast}}\E\vert\Delta_0\Delta_i\vert\\
&\leq\frac{\kappa\,\E(\Delta_0^2)}{(\vert\Lambda_n\vert b_n)^{1/2}}+\left(\frac{\kappa}{(\vert\Lambda_n\vert b_n)^{1/2}}+2\right)\sum_{i\in W_n^{\ast}}\E(\vert\Delta_0\Delta_i\vert)\\
&=O\left(\frac{1}{(\vert\Lambda_n\vert b_n)^{1/2}}+\left(\frac{1}{(\vert\Lambda_n\vert b_n)^{1/2}}+2\right)\left(N_n^db_n+o(1)\right)\right)\qquad\textrm{(by Lemma $\ref{esperance_delta0_deltai}$)}\\
&=o(1)\qquad\qquad\quad\textrm{(by (\ref{Nn}) and Assumption $\textbf{(A2)}$)}.
\end{align*}
In order to obtain ($\ref{equation1}$) it remains to control
$$
\textrm{F}_{0}=\E \left(\sum_{k=1}^{\vert\Lambda_n\vert}h_{k-1,k+1}^{''}\left(\frac{\Gamma_{\varphi(k)}^2}{2}
+\Gamma_{\varphi(k)}\left(S_{\varphi(k-1)}(\Gamma)-S_{\varphi(k)}^{(n)}(\Gamma)\right)-\frac{1}{2\vert\Lambda_n\vert}\right)\right).
$$
Applying Lemma $\ref{esperance_delta0_deltai}$, we have
\begin{align*}
\textrm{F}_0&\leq\left\vert\E\left(\frac{1}{\vert\Lambda_n\vert}\sum_{k=1}^{\vert\Lambda_n\vert}h_{k-1,k+1}^{''}(\Delta_{\varphi(k)}^2-\E(\Delta_0^2))\right)\right\vert
+\vert1-\E(\Delta_0^2)\vert+2\sum_{j\in V_0^1\cap W_n}\E\vert\Delta_0\Delta_{j}\vert\\
&\leq\left\vert\E\left(\frac{1}{\vert\Lambda_n\vert}\sum_{k=1}^{\vert\Lambda_n\vert}h_{k-1,k+1}^{''}(\Delta_{\varphi(k)}^2-\E(\Delta_0^2))\right)\right\vert
+O(N_n^db_n)+o(1).
\end{align*}
Since $N_n^db_n\to0$, it suffices to prove that
$$
\textrm{F}_1=\left\vert\E\left(\frac{1}{\vert\Lambda_n\vert}\sum_{k=1}^{\vert\Lambda_n\vert}h_{k-1,k+1}^{''}(\Delta_{\varphi(k)}^2-\E(\Delta_0^2))\right)\right\vert
$$
goes to zero as $n$ goes to infinity. In fact, keeping in mind $W_n=\{-N_n+1,...,N_n-1\}^d$ and $W_n^{\ast}=W_n\backslash\{0\}$, we have 
$$
\textrm{F}_1\leq\frac{1}{\vert\Lambda_n\vert}\sum_{k=1}^{\vert\Lambda_n\vert}\left(\textrm{L}_{1,k}+\textrm{L}_{2,k}\right)\\
$$
where $\textrm{L}_{1,k}=\left\vert\E\left(h_{k-1,k+1}^{''(n)}\left(\Delta_{\varphi(k)}^2-\E\left(\Delta_0^2\right)\right)\right)\right\vert=0$ 
since $h_{k-1,k+1}^{''(n)}$ is $\F_{V_{\varphi(k)}^{m_{n,\infty}}}$-measurable and
\begin{align*}
\textrm{L}_{2,k}&=\left\vert\E\left(\left(h_{k-1,k+1}^{''}-h_{k-1,k+1}^{''(n)}\right)
\left(\Delta_{\varphi(k)}^2-\E\left(\Delta_0^2\right)\right)\right)\right\vert\\
&\leq\E\left(\left(2\wedge\sum_{i\in W_n}\frac{\vert\Delta_i\vert}{\vert\Lambda_n\vert^{1/2}}\right)\Delta_0^2\right)\\
&\leq \kappa\left(\frac{\E(\Delta_0^2)}{(\vert\Lambda_n\vert b_n)^{1/2}}+\frac{\sum_{i\in W_n^{\ast}}\E\vert\Delta_0\Delta_i\vert}{(\vert\Lambda_n\vert b_n)^{1/2}}\right)\quad\textrm{(since $\vert\Delta_0\vert\leq\frac{\kappa}{\sqrt{b_n}}$ a.s.) }\\
&=O\left(\frac{1}{(\vert\Lambda_n\vert b_n)^{1/2}}+\frac{N_n^db_n+o(1)}{(\vert\Lambda_n\vert b_n)^{1/2}}\right)\qquad\textrm{(by Lemma $\ref{esperance_delta0_deltai}$)}\\
&=o(1)\qquad\qquad\textrm{(by (\ref{Nn}) and Assumption $\textbf{(A2)}$)}.
\end{align*} 
The proof of Theorem $\ref{convergence-loi}$ is complete.\\
\\
\textbf{Acknowledgments}. The author thanks an anonymous referee for his$\backslash$her careful reading and constructive comments.

\bibliographystyle{plain}
\bibliography{xbib}
\end{document}